\documentclass[11pt]{article}

\usepackage{cite}

\begin{document}

\title{The history of V.~A.~Rokhlin's ergodic seminar (1960--1970)}
\author{A.~M.~Vershik\footnote{St.Petersburg Department of Steklov Institute of Mathematics, St.Petersburg, Russia,
email: {\tt avershik@pdmi.ras.ru}}}

\date{}

\maketitle

\begin{abstract}{The paper tells about the main features and events of the ergodic seminar organized and
headed by V.~A.~Rokhlin at the Leningrad State  University. The seminar was active in {\rm1960--1970.}}
\end{abstract}

\section*{1. V.~A.~Rokhlin's moving to Leningrad}

Vladimir Abramovich Rokhlin, with his wife Anna Alexandrovna Gurevich, his son Volodya (born in 1952),
and his daughter Lisa (born in 1955),  moved to Leningrad in the late summer of 1960. Before, they lived
in Kolomna near Moscow, where V.A.\ worked at a pedagogical institute. In 1959, the rector of the
Leningrad State University A.~D.~Alexandrov offered V.A.\  a~position of professor at the Chair
of Geometry, of which he was the head. The corresponding advice Alexandrov received  from a friend and
colleague N.~V.~Efimov, a well-known Moscow geometer, whose family were old friends of  the Rokhlins:
Efimov's wife was a friend of A.A.\ from the 1930s, when they both were L.~S.~Pontryagin's graduate
students in Voronezh.

Immediately upon arrival in Leningrad, in September 1960, V.A.\ announced his research  seminar on
the theory of dynamical systems (``ergodic seminar''). Besides, after a while he started  to teach an
elective course in combinatorial topology. Neither the first area of mathematics (the theory of
dynamical systems), nor the second one (modern topology) was represented at the university before.
Another seminar, on topology, was announced by V.A.\ somewhat later, since it was necessary first to
train at least the first group of participants. This seminar was started in 1961, became well known, and
existed almost until V.A.'s death in 1984. It deserves a separate memoir.

\section*{2. The prehistory}

Here I want to present the history of the ergodic seminar. Its rapid creation and immediate activity
are explained, first, by the fact that it coincided with an impressing progress in the theory of
dynamical systems caused by the appearance of Kolmogorov's entropy and other events, succeeded by an
enormous interest to the subject. But there was also another circumstance: the ground for the seminar
was well prepared. V.A.'s arrival in Leningrad was preceded by one year by the arrival of his first
graduate student  in dynamical systems Leonid Mikhaylovich Abramov, whose wife was from Leningrad. In
1955, Abramov graduated from the Dnepropetrovsk University with a degree in approximation theory and was
sent to work at the Arkhangelsk Pedagogical Institute, where V.A.\ held a position from 1952, when he
was forced to leave Moscow after he had been effectively expelled from the Steklov Institute of
Mathematics by its director I.~M.~Vinogradov.

L.~M.~Abramov, on V.A.'s advice, began to study ergodic theory, and in the mid-1950s enrolled in the
postgraduate program of the Moscow State University under V.A.'s supervision  (by agreement with
A.~N.~Kolmogorov). L.M.'s thesis became the first PhD thesis on entropy theory (defended at the Moscow
State University in 1960); he himself was the first V.A.'s graduate student. Many knew about
V.A.'s plans of moving to Leningrad, some (me included) knew from Moscow friends about L.M.'s arrival.
Several people, including O.~A.~Ladyzhenskaya, who was friends with V.A.\ from the time of her
university studies in Moscow, helped L.M.\ to solve the problem, extremely difficult at the time, of
getting a job at the Leningrad State University. Coincidentally, at this very time L.~V.~Kantorovich
organized a new department of mathematical economics at the Faculty of Economics, and this made it
possible to get L.M.\ a job at the Chair of Mathematical Economics, where he worked till the end of
his life. The times were like that.

In Leningrad, at G.~P.~Akilov's  (1921--1986) seminar, there was interest in measure theory on
locally compact groups and especially in measure theory on linear topological spaces (the influence
of I.~M.~Gelfand and D.~A.~Raikov). The core members  of the seminar were B.~M.~Makarov, V.~P.~Khavin,
V.~N.~Sudakov, and myself. Some of us had already heard about Kolmogorov's entropy. In the fall of 1956,
the Second All-Union Conference on Functional Analysis was held in Odessa.  Like the first one (which took place  in Moscow in January 1956), it attracted many well-known specialists. In particular, S.~V.~Fomin gave a talk about  the famous Kolmogorov's paper on entropy ({\it Dokl. Akad. Nauk SSSR}, 1958), which had just appeared. At the next, third, conference on functional analysis, which was held in Baku in the spring of 1960, V.A.\ gave a talk about the  progress in the theory of dynamical systems related to the notion of entropy. Together with other graduate students from Leningrad (B.~M.~Makarov,  V.~N.~Sudakov), I participated in both conferences. We were very impressed by these talks, and they aroused a great interest in the subject.

An interesting episode occurred at the Odessa conference. After S.~V.~Fo\-min's talk about Kolmogorov's entropy, an unknown graduate student of the Odessa University D.~Z.~Arov unexpectedly rose from his seat to claim that A.~N.~Kolmogorov's results, presented in the talk, on applications of Shannon's entropy to Bernoulli schemes were contained in his master thesis, of which he informed A.N.\ in a letter. A.N.'s second paper on entropy ({\it Dokl. Akad. Nauk SSSR}, 1959) does indeed mention D.~Z.~Arov's letter. But what exactly had been done by D.~Z.~Arov and how this is related to the Kolmogorov--Sinai entropy became clear only quite recently. The part of his master thesis related to information theory was recently published in our series of {\it Zapiski Nauchnykh Seminarov POMI} (Vol.~436, 2015) with    commentaries written, at my request, by B.~M.~Gurevich. Thus, D.~Z.~Arov may be safely called one of the pioneers of entropy ergodic theory.

Later, in 1962, D.~Z.~Arov came to Leningrad in order to work with V.A.\ and took an active part in the ergodic seminar, but then he switched to a field close to his teacher M.~G.~Krein.

In the 1959--1960 academic year, we organized a home student seminar on measure theory and functional analysis, where we studied books by A.~Weil, P.~Halmos, and some papers on integration theory and measure theory.

It is appropriate to recall V.A.'s visit to Leningrad in 1957 to give a talk at G.~M.~Fikh\-tengolts's and L.~V.~Kantorovich's seminar. I was present at this meeting and remember the talk: it was one of the last meetings of the seminar, soon it ceased to exist because of L.V.'s moving to Novosibirsk. The talk was devoted to the recent V.A.'s paper in {\it Uspekhi Mat. Nauk} on metric invariants of measurable functions based on his theory of measurable partitions. I remembered this talk, and many years later used and generalized it in several directions. During the same visit (or another one close in time), V.A.\ gave a lecture at the Leningrad Department of Steklov Institute of Mathematics on recent progress in algebraic and combinatorial topology. A rotaprint edition of  the text of the lecture was immediately issued, and for a long time it attracted the attention of young mathematicians.

\section*{3. The early days of the seminar and its participants}

Upon arrival, V.A.\ asked L.~M.~Abramov, who by September 1960 already had good contacts with Leningrad mathematicians, to put him in touch with those wishing to take part in an ergodic seminar. The first organizational meeting, as far as I remember, took place not later than October 1960 at the Faculty of Mathematics and Mechanics. L.M.\ became the secretary of the seminar.

The fact that V.A.\ was now a professor at the  Faculty of Mathematics and Mechanics was known to many; the interest to him and to the new seminar was great. That is why, the first meeting attracted many undergraduates, graduate students, and professors. V.A.\ proposed to begin
 with the study of    measure theory in his understanding, the theory of Kolmogorov's entropy, as well as a series of classical works on dynamics: spectral theory, automorphisms of compact groups, geodesic flows, relation to stationary random processes, etc. All these subjects were a complete novelty for most of the audience. Now I remember only those who took part in the seminar for a sufficiently long time in    its early years. Among them, I name only a few:
L.~M.~Abramov, I.~A.~Ibragimov, V.~N.~Sudakov, B.~M.~Makarov, A.~M.~Kagan, R.~A.~Zaidman. Later, during the 1960s, some young mathematicians, who had graduated from the University or the postgraduate program,
 joined in: R.~M.~Belinskaya, S.~A.~Yuzvinsky, M.~I.~Gordin, A.~A.~Lodkin, A.~N.~Livshits, and others. Of course, there were those who had just flickered and disappeared, there were, especially at first, curiousity seekers; the corresponding list is long and difficult to recover.

The first meeting of the seminar was opened by V.A., who told about his plans for the seminar and suggested the topics of the first talks (see below). These topics were rather quickly taken up, and soon regular meetings began. The first invited speakers included two leading Leningrad mathematicians who worked, in particular, on problems close to dynamics and  with whom V.A.\ maintained close contact: D.~K.~Faddeev, who spoke about the classification of automorphisms of lattices, and Yu.~V.~Linnik, who gave a general outline of his recent paper on the equidistribution of integer points on manifolds. I have kept a copy of the first paper (later  transformed into the book {\it Ergodic Method and $L$-Functions}), which was presented by Yu.V.\  to V.A. and passed to me by V.A.'s widow,  bearing the inscription: ``This is my best paper.''

An important feature of the seminar was the abundance of guests, mainly from Moscow, where Ya.~G.~Sinai's and V.~M.~Alekseev's seminar on dynamical systems was very active at the time. This seminar followed A.~N.~Kolmogorov's initiative of the late 1950s to revive the theory of dynamical systems in Moscow and, in particular, to assert entropy theory, which was emerging after his pioneering works,  as one of its topics. At A.N.'s probability seminar,  this topic received some (little) attention, and V.A., which at the time worked near Moscow, took part in it. He immediately appreciated A.N.'s idea and became one of the main developers of this theory. Recall that in his second note on entropy, A.N.\ emphasizes the role of the theory of measurable partitions created by V.A.\ and in a footnote cites the counterexample suggested by V.A.\ to one of the statements of the first note. V.A.\ repeatedly mentioned the very great impression which Kolmogorov's work on entropy made on him. To this work he attributed his own renewed interest  to the theory of dynamical systems, the subject of  his famous papers of the late 1940s, which made him (along with the subsequent works on entropy theory) a leading expert in this field.

All of these things largely predetermined also the topics of the Leningrad ergodic seminar. And they explain the fact that its guests and speakers included all active Moscow ``dynamists'' of the 1960s: V.~I.~Arnold, Ya.~G.~Sinai, D.~V.~Anosov, V.~M.~Alekseev, and somewhat later their pupils B.~M.~Gurevich, V.~I.~Oseledets, A.~B.~Katok, A.~M.~Stepin, G.~A.~Margulis, and others. There were also guests from other cities: Odessa, Tashkent, Tbilisi, Novosibirsk, Gorky, etc.

The seminar existed up to the late 1960s. Occasional meetings (say, on the occasion of visits or special events) took place even later. I remember that I spoke about my Doctor of Sciences dissertation in 1971. In 1965, I started my own seminar on smooth dynamics, whose participants included some of those I have already mentioned, and also M.~L.~Gromov, V.~L.~Eidlin, S.~M.~Belinsky, and others. The seminar existed for one and a half year. In 1968, I organized a wider seminar on measure theory in linear spaces, operator algebras, representation theory, and dynamics, which exists to this day.

\section*{4. Topics of the seminar}

 The main topic of the ergodic seminar was, of course, ergodic theory rather than the entire dynamics. V.A.\ emphasized that he advocated distinguishing between mathematical structures and preferred not to mix, say, smooth dynamics with metric or topological dynamics. At the same time, he knew classical dynamics (both smooth and number-theoretic) well, but for him they existed separately, and he gave precedence to metric (measure-theoretic) dynamics;  I dare say, for aesthetic reasons. Moreover, he believed, for example, that it is the methods of the theory of measure-preserving transformations that would eventually allow one to solve many problems of both analytic number theory and smooth dynamics.

He particularly emphasized the role of dynamical systems of algebraic origin as natural examples. Topological dynamics, as compared with metric dynamics, was for him a less natural area of mathematics. More exactly, he put emphasis on the existence of an invariant measure with respect to a~group of transformations. Of course, such a measure does not always exist, but then this is rather a flaw in the problem formulation, and it is the invariance of a measure with respect to some group of transformations that ensures the meaningfulness of studying this measure. On the other hand, V.A.\ believed that measure theory itself is not an especially important field, in particular, thought that the theory of nonseparable measure spaces is a~nonessential and somewhat pathological object. He believed that Lebesgue spaces, introduced and studied by him, cover the needs of analysis and dynamics in measure spaces. The widespread clich\'e ``measure space'' without further details invited his criticism. He preferred the more accurate term ``the theory of transformations with an invariant measure'' to the vague  term
 ``ergodic theory.'' I plan to express my views on the role of V.A.'s seminal paper ``On the fundamental ideas of measure theory'' and on the modern understanding of what is the category of measure spaces elsewhere.

Concluding the general description of the seminar, I want to list several principles which, as I think, guided V.A.\ in mathematics in general and, in particular, in his style of
conducting a research mathematical seminar. First of all, he insisted that all statements must be clear-cut and, if possible,
structurally (or categorically) pure. Further, he insisted
on the invariance (functoriality) of statements and, again if possible, proofs. He rejected arguments that went beyond the framework of a given structure or category, and this occasionally caused discussions, which usually ended in his favor. This attitude could be regarded as an extreme kind of Bourbakism, but V.~I.~Arnold, a well-known enemy of Bourbaki, eagerly supported V.A.'s recommendations and followed them.

Besides, V.A.\ was able to appreciate the beauty of clear-cut statements and the originality of arguments, and this feature is in a sense opposite to any dogmatism. Anyway, the school of V.A.'s seminar had much to offer to those who went through it and, possibly, did not subscribe to all his principles.

\section*{5. On the style of conducting the seminar and the first meetings}

Unfortunately, the agenda of all  seminar  meetings for almost ten years of its existence has not survived. Now it seems astonishing: thanks to modern electronic techniques, we easily save anything and everything. But at the time, it was necessary to keep records, follow the sequence of talks, which only few were capable of. As a result, we (and the history of science) have lost very much. I wrote about this in relation to the famous I.~M.~Gelfand's seminar (see M.~A.~Shubin's notes of talks at the seminar available in the internet); but  this applies as well to other seminars.

The above-mentioned numerous guests of the seminar spoke about their published papers; it is hardly worthwhile (and anyway impossible) to list them or reproduce the corresponding discussions. Of course, V.A.'s opinion  on a result or a talk, with which he usually concluded the meeting, as well as his judgements and comments, were always of interest both to the audience and speakers. He always expressed these opinions with authority and reserve  no matter whether the talk was received with enthusiasm or  with doubt. It is worth mentioning that V.A.\ disapproved of other styles of conducting a seminar, when the leader, albeit a distinguished scholar, reacts to talks with unrestrained emotions.  The culture of mathematical seminars in Soviet time was being worked out gradually, following, apparently, European and pre-revolutionary traditions. V.A.\ followed the best of them.

It remains to extract from my memory what it has preserved. Unsurprisingly, I remember mainly seminars somehow related to subjects of interest to myself, and some of my own talks. To begin with, I present only a partial list of talks of the first period about published results given mainly by the participants of the seminar.

During the initial period, in 1960--1961, the talks were about   famous papers on ergodic theory and measure theory: about the Lebesgue space (V.~N.~Sudakov), the definitions of entropy (I.~A.~Ibragimov), spectral theory (L.~M.~Abramov, quasidiscrete spectrum; B.~M.~Makarov, the spectrum of Gaussian systems), the martingale theorem and ergodic theorem (A.~M.~Kagan), etc.

I will write in slightly more detail about my own talks of the initial period. First, V.A.\ asked me to  tell about his little-known paper on unitary rings, which he valued. It contained an attempt, in the spirit of Gelfand--Naimark and Moscow functional analysis, to present measure theory and the notions of ergodic theory in terms of the space $L^2$ of square-integrable functions. For this,  the Hilbert structure of the space of functions had to be equipped with some additional structure, and V.A.\ choose a structure which he called a unitary ring (later, it became known as a unitary algebra). In other words, one considers an axiomatization of an unbounded (i.e., defined not for all pairs of elements) multiplication in~$L^2$ with obvious constraints. Indeed, this defines a functor between the corresponding categories (the category of Lebesgue spaces and the category of unitary rings; the terms did not existed at the time). To an automorphism with an invariant measure there corresponds a multiplicative unitary operator. This is a natural measure-theoretic counterpart of Gelfand's theory of normed rings. However, another category  has become more popular here, whose objects are $L^2$~spaces with a cone, the cone of nonnegative functions. This functor proved to be more convenient, for example, when dealing with operators: nonnegative automorphisms or (later) polymorphisms and Markov operators. But my talk contained an answer to a specific question indirectly posed by V.A.: how can one define a structure of a unitary ring on an infinite-dimensional space with a Gaussian measure?

At the time, V.~N.~Sudakov and myself were interested mainly in geometric aspects of
Gaussian measures  (quasiinvariant measures, extension of weak distributions, etc.), partly under the influence of the Moscow school (I.~M.~Gelfand, D.~A.~Raikov). In particular, I studied N.~Wiener's book {\it Nonlinear Problems in Random Theory} (1958). But the interest to the dynamics of linear automorphisms of Gaussian measures, or, as we said back then, ``normal dynamical systems,'' was initiated by A.~N.~Kolmogorov and picked up by V.A., S.~V.~Fomin, and others. This applied to spectral theory, entropy, etc. A.N.'s remarkable idea, implemented by his pupil I.~V.~Girsanov, led to constructing an automorphism with simple singular spectrum (in the complement to the constants). V.A.\ popularized this field, those who worked in it included
L.~M.~Abramov    and later Ya.~G.~Sinai and many others,   including myself. In my talk, I described multiplication formulas for generalized Hermite polynomials which defined the structure of a unitary ring. Then, and somewhat later, I understood that these formulas are equivalent to formulas from the It\^o--Wiener stochastic analysis, or Wick's regularization formulas (in physics), etc.

I remember that V.A.\ liked these results very much. Later, they were included in my PhD thesis. Back then, I was a graduate student, formally under the supervision of G.~P.~Akilov, a pupil and coauthor of L.~V.~Kantorovich, to whom I owe much, but actually my thesis was inspired by V.A. Among other theorems proved in the thesis, it announced a result on the isomorphism of Gaussian automorphisms with absolutely continuous spectrum, which was proved only for some special spectra. In full generality it was proved somewhat later by D.~Ornstein. Afterwards, I repeatedly spoke about  this   at the seminar.

I remember some of my other talks, initiated by V.A., about famous papers, for example, my talk about  the paper  ``Geodesic flows on manifolds of constant negative curvature'' by
I.~M.~Gelfand and S.~V.~Fomin ({\it Uspekhi Mat. Nauk}, Vol.~7, No.~1(47), 1952) and the corresponding discussion, as well as my talk about  S.~Smale's remarkable paper on the orthogonal group as a retract of the group of homeomorphisms.

\section*{6. Talks of    seminar participants about their own results}

Almost all participants of the ergodic seminar worked, apart from ergodic theory, also in other fields of mathematics, and even reported the obtained results at the seminar. I remember V.~N.~Sudakov's talks about measures in functional spaces and my own talks about the axiomatics of measures in linear spaces from the standpoint of Lebesgue spaces, etc. I had already begun to work in representation theory, algebraic and combinatorial asymptotics, always bearing in mind connections and analogies with dynamics. But, speaking about not educational or ``exterior'' talks, but research talks on the topic of the seminar, irrespective of the time when they were given, I must first of all mention, at least very briefly, the most important results of the 1960s--1970s obtained by the permanent participants of the seminar in ergodic theory proper.

1. V.A.: a proof of the existence of a countable generator for every aperiodic endomorphism; later, he proved the  existence of a generator with finite entropy for endomorphisms with finite entropy. The ultimate result on the existence of a finite generator for such endomorphisms with an estimate was proved by W.~Krieger and, independently, by my graduate student A.~N.~Livshits. But V.A.'s result was the first in this direction.

2. V.A.\ and L.~M.~Abramov: a formula for the entropy of skew products using the new notion of the mixed entropy of fibers.

3. V.A.\ and Ya.~G.~Sinai: a proof of the fundamental theorem on the coincidence of the classes of Kolmogorov automorphisms and  of automorphisms with totally positive entropy; this is one of the fundamental facts of ergodic theory. Later, it was extended by D.~Ornstein and B.~Weiss to actions of amenable groups.

The next two permanent participants of the seminar from the mid-1960s deserve special mention. During the years of V.A.'s teaching at the Leningrad State University, he had many PhD and master students in topology, but, apart from his very first PhD student L.~M.~Abramov, only two PhD students in ergodic theory. The first of them was S.~A.~Yuzvinsky, who brilliantly completed his evening studies at the Faculty of Mathematics and Mechanics of the Leningrad State University (earlier, he had graduated from the Polytechnic Institute), but nevertheless, in spite of V.A.'s recommendation,  was not admitted to the postgraduate program at the University. He enrolled in such a program at the Hertzen Institute, where V.A.\ worked part-time for a while, and under V.A.'s supervision prepared and defended his PhD thesis (in 1966).

4. S.~A.~Yuzvinsky's results on the entropy theory of automorphisms of compact groups. This was one of V.A.'s favorite topics. He himself and many his pupils and colleagues (Ya.~G.~Sinai and his pupils, D.~Z.~Arov, L.~M.~Abramov) studied entropy formulas for group automorphisms. S.\,A.'s results and subsequent work by K.~Schmidt and others (on the Mahler measure) in a sense settled this problem. Another specific result obtained by S.\,A., on the genericity  of  automorphisms with simple continuous spectrum, also belongs to the circle of problems traditional for approximation theory in ergodic theory.

5. The second V.A.'s PhD student in ergodic theory was R.~M.~Belinskaya (Ekhi\-lev\-s\-kaya, 1938--2011).  Her main topic was time changes in automorphisms and related problems of ergodic theory. The results of her thesis (1970) are widely used in the literature. The main  problem suggested to her by V.A.\ will be described below. But first it is worth mentioning that V.A., at our request,  made several surveys of open problems in ergodic theory, first time around 1962, and then in 1965 (mostly about entropy-related problems).

V.A.\ mentioned various problems, but all of them were related to his own work. In the first talk, he promoted the study of endomorphisms and  suggested to classify finite monotone sequences of measurable partitions as their invariants. It is appropriate to mention, in particular, O.~V.~Guseva, a participant  of the seminar   during its early years,  who worked mainly on the theory of partial differential equations, but attended the ergodic seminar, took up this subject, and, generalizing the classification of measurable partitions given by V.A., obtained a complete classification in the case of finite sequences.

Another problem was especially important: V.A.\ asked when two ergodic automorphisms have isomorphic orbit partitions. He said that  this question appeared in his studies in the 1940s, but it was not known even whether this is true for two nonisomorphic rotations of the circle. V.A.\ himself believed that in this case the answer should be positive. And his bold conjecture was that if two automorphisms have different entropies, then their orbit partitions are nonisomorphic. It is this problem that R.~M.~Belinskaya systematically tried to solve for a long time, discovering interesting facts about time changes along the way. In particular, she observed that the orbit partition of an automorphism is the intersection of a decreasing dyadic sequence of partitions. But the main problem remained open. At the end of 1966, R.M.\ and S.A.\ informed me of this fact, and soon I proved that all ergodic dyadic sequences are lacunarily isomorphic, and hence their intersections are metrically isomorphic. Thus, V.A.'s conjecture turned out to be false. V.A.\ was very glad that the answer was finally obtained. The results were published in {\it Functional Analysis and Its Applications} in 1968. But while the paper was in preparation, I discovered, selecting material for our survey, joint with S.A., on dynamical systems and operator algebras for {\it Itogi Nauki}, that as early as in 1963, the American mathematician H.~Dye proved this fact in terms of the theory of operator algebras etc. Dye mentioned that this problem had  been implicitly posed by J.~von Neumann, in a paper from the famous series of joint papers with F.~Murray, as the isomorphism problem for hyperfinite factors of type~${\rm II}_1$. Von Neumann believed that the answer is most probably positive, but that it would require  a detailed analysis of ergodic automorphisms. It is this problem that Dye solved. It is worth mentioning that V.A.\ knew about the  isomorphism problem for factors, but did not know about its relation to his question. Dye's theorem is a classical fact of ergodic theory and is much cited in the literature. My proof, relying on the lacunary isomorphism theorem, is much simpler than the original one. Of course, these results were reported both at the ergodic and Moscow seminars.

At first, it seemed to me that a strengthening of the method of proving the lacunary isomorphism theorem would allow one to prove the isomorphism of any two ergodic dyadic sequences. In my paper on lacunary isomorphism published in {\it Functional Analysis}, I even made the corresponding remark. However, very soon I became aware of the difficulty of the problem, and it took me a year before I managed, on the one hand, to find a ``standardness criterion,'' i.e., a necessary and sufficient approximative condition for a given dyadic sequence to be isomorphic to a simplest (or standard) sequence, and, on the other hand, to present a continuum of pairwise nonisomorphic ergodic dyadic sequences (they are now called dyadic filtrations). The corresponding invariant was a new metric invariant, the entropy of a~filtration. Following the natural logic of the subject, I defined the so-called ``main invariant'' of a~filtration, and somewhat later, the corresponding invariant of an action, the ``scale'' of an automorphism. Soon after my results were reported at the Moscow ergodic seminar, A.~M.~Stepin proved that the Kolmogorov entropy of the action of the group~$\sigma_1^{\infty} Z_2$ is also an invariant of a dyadic sequence. The entropy of a filtration introduced by me coincides in this case with the entropy of the action of the group. Curiously enough, V.A.'s intuition has been eventually almost vindicated: the entropy of an action is an invariant close to an orbit invariant, but still not exactly an orbit invariant. The main results were included into my Doctor of Science dissertation ``Approximation in measure theory.''  As far as I remember, it was reported at probably the last meeting of the ergodic seminar (around~1971).

Later, this circle of problems (orbit theory, the theory of filtrations, etc.) grew into a large research area and found applications in a vast variety of fields: the theory of random processes, asymptotic combinatorics, and, of course, ergodic theory. The theory of filtrations must provide the foundation for a combinatorial (approximation) approach to dynamics. At the moment, it is far from completed.

A few words about other seminar participants.

I have already mentioned the young mathematicians M.~I.~Gordin,\break A.~A.~Lodkin, A.~N.~Livshits, and others, who took an interest in ergodic theory. Among my graduate students of that time working in this field, A.~N.~Livshits (1950--2008) stood out. He was a highly gifted and early matured mathematician. He had caught only the very end of the ergodic seminar, and participated mainly in my later seminar, to which he came as a~first-year student.  Very quickly A.~N.~Livshits mastered the theory of Anosov systems and, while a third-year undegraduate, proved his famous theorem on the cohomology of hyperbolic automorphisms and flows. Then he worked much on encoding and isomorphism, independently proved Krieger's theorem (see above), and later we jointly established  a direct relation between adic transformations and substitutions.

By the time of the new flourishing of ergodic theory connected with  the name of D.~Ornstein, who solved the isomorphism problem for Bernoulli automorphisms and discovered non-Kolmogorov Bernoulli systems, the seminar had already slowed its activity. These papers were not discussed at the seminar. But it is appropriate to mention that one of the seminar participants  R.~A.~Zaidman, prior to Ornstein's works, claimed to have proved that the entropy is  a complete invariant of Bernoulli automorphisms. This claim was  widely known. However, repeated attempts of seminar participants to understand R.~A.~Zaidman's arguments failed. It is worth mentioning that he is a very talented and interesting mathematician, and his ideas on encoding were later picked up by A.~N.~Livshits, who greatly developed the theory of codes.

The works of other my graduate students in ergodic theory were done later.

In 1969, V.A.\ offered me to conduct, together with him, a small ergodic seminar for beginners (first and second year students), and I gladly accepted. V.A.\ only outlined the program of the seminar, but almost did not attend its meetings. The seminar was educational and was active for one term. Its participants and speakers included the students A.~A.~Suslin, V.~V.~Rokhlin~(Jr.), A.~G.~Reiman, V.~M.~Kharlamov, and many others, who later became well-known mathematicians. By the way, this generation (graduated in 1972) was one of the strongest at the Faculty of Mathematics and Mechanics in those years (A.~N.~Livshits, B.~S.~Tsirelson, A.~A.~Suslin, V.~M.~Kharlamov, and others).

\section*{7. Schools, conferences, the conclusion}

Speaking about the ergodic seminar in the 1960s, one cannot fail to mention the ergodic school (August 1965) at Khumsan (Uzbekistan). It took place during the heyday of ergodic theory itself and relatively liberal times in the USSR. The venue and conveniences were excellent.  The organizers were mainly from Moscow (Ya.~G.~Sinai and his pupils) and Uzbekistan (T.~A.~Sarymsakov and others). The school attracted an exceptional group of participants and speakers, including  a very wide participation by experts in other fields of mathematics. Here is an incomplete list: V.~A.~Rokhlin, O.~A.~Ladyzhenskaya,   T.~A.~Sarymsakov, N.~P.~Romanov, Yu.~M.~Smirnov, V.~P.~Maslov, A.~Ya.~Povzner, I.~I.~Pyatetsky-Shapiro, Ya.~G.~Sinai, A.~A.~Kirillov, many of them with families. Besides, a whole team of young talents arrived, mainly from Moscow: I.~N.~Bernstein, M.~L.~Gromov, A.~B.~Katok, A.~G.~Kushnirenko, L.~G.~Makar-Limanov, G.~A.~Margulis, A.~M.~Stepin, and others.

In conclusion, I should mention several more facts from the history of the ergodic seminar. In August 1969, we organized a one-day informal conference dedicated to V.A.'s 50th anniversary. Numerous participants listened with interest to talks given by V.~I.~Arnold, S.~P.~Novikov, Ya.~G.~Sinai, D.~V.~Fuks, M.~L.~Gromov, and A.~M.~Vershik. Ten years later, on December~10,~1979, a similar one-day conference was organized  in honor of V.A.'s 60th anniversary, with the same speakers, except that instead of  M.~L.~Gromov's talk there was V.~M.~Buchstaber's talk. The day was chosen by chance, the organizer forgot that in the free world calendar  December 10 is the Human Rights Day, and that in the USSR every mention of this day,  let alone the organization of any meetings on it, is punishable. That is why, as I was later told by V.~A.~Zalgaller,  the Party bosses of the faculty were terribly frightened and ran to the District Party Committee (!) for consultations. But eventually it was decided that  it was safer for all involved to pretend that there was no conference at all.

Soon after that, the faculty administration suggested that V.A.\ retire, under a pretext that he declared that for health reasons he cannot teach at the new building of the faculty in New Peterhof, but can give seminars in the city (for example, at the Leningrad Department of Steklov Institute of Mathematics). After much efforts and letters to  the university and faculty administrations, he was allowed to continue working at the Faculty of Mathematics and Mechanics, but only for one year, as ordered by the then rector of the Leningrad State University (V.~B.~Aleskovsky).

In V.A.'s notes, which I found after his death (December 3, 1984), the last manuscript concerning ergodic theory was dated 1968 (classification of measurable partitions); the text was unfinished. Later, V.A.\ worked only in topology and real algebraic geometry, but followed the progress in ergodic theory and attended the seminar up to the very last meetings.

\end{document}